# Capacities and Hessians in a class of $m$–subharmonic functions.

A. S. Sadullaev and B. I. Abdullaev

**0.** For a twice differentiable function $u \in C^2(D)$ in a domain $D \subset \mathbb{C}^n$ differential operator $(dd^c u)^m \wedge \beta^{n-m}$ called *hessian* of function $u$, where $d = \partial + \bar{\partial}$, $d^c = \frac{\partial - \bar{\partial}}{4i}$ and $\beta = dd^c |z|^2 -$ standard form of volume in $\mathbb{C}^n$. The notation of the operator justified that, if $\lambda(u) = (\lambda_1(u), ..., \lambda_n(u)) -$ vector of eigenvalues of hermitian matrix $(u_{j\bar{k}})$ of quadratic form $dd^c u = \frac{i}{2}\sum_{j,k} u_{j\bar{k}} dz_j \wedge d\bar{z}_k$, then

$$(dd^c u)^m \wedge \beta^{n-m} = m!(n-m)! H_m(u)\beta^n \ , \qquad (1)$$

where $H_m(u) = \sum_{1 \le j_1 < ... < j_m \le n} \lambda_{j_1} ... \lambda_{j_m}$ – Hessian of vector $\lambda(u) \in \mathbb{R}^n$.

The aim of this paper is to study $m-$subharmonic functions connected with operator (1) and an equation

$$(dd^c u)^m \wedge \beta^{n-m} = f(z)\beta^n \ , \qquad (2)$$

also, construction of potential theory on their basis. At $m = 1$ the equation (2) gives a Poisson equation and at $m = n$ it gives a Monge –Ampere equation; which good developed and constitute of fundamentals of classical and complex potential theory. In general case $1 \le m \le n$ the equation (2) called complex equation of Hessian. This equation and properties of their solutions was studied systematically in the past ten years. Here we bring a reference only to some works, which has directly relations on this paper [5,8,9,12-15] , especially , Z. Blocki [5] and Dinev S., Kolodziej S. [8] from which we take main symbols and methods of studying of Hessians.

**1. Definition 1.** *Twice differentiable function* $u \in C^2(D)$, $D \subset \mathbb{C}^n$ *called* $m-$*subharmonic* $(m-sh)$ *in* $D$ $(1 \le m \le n)$, *if*

$$(dd^c u)^k \wedge \beta^{n-k} \ge 0 \ , \ \forall \ k = 1,2,...,m \ . \qquad (3)$$

We have following statement:

$$dd^c u_1 \wedge dd^c u_2 \wedge ... \wedge dd^c u_m \wedge \beta^{n-m} \ge 0 \ \ \forall \ u_1, u_2, ..., u_m \in m - sh(D) \cap C^2(D) \ . \qquad (4)$$

This statement has a dual character: if $u_1$ twice differentiable and satisfy (4) for all $u_2,...,u_m \in m-sh(D)\cap C^2(D)$, then it is $m-sh$. This condition allow us define $m-sh$ functions in the class of $L^1_{loc}$ functions.

**Definition 2.** *A function $u \in L^1_{loc}(D)$ called $m-sh$ in $D\subset\mathbb{C}^n$, if it is upper semicontinuous and for any twice differentiable $m-sh$ functions $v_1,...,v_{m-1}$ a current $dd^cu\wedge dd^cv_1\wedge...\wedge dd^cv_{m-1}\wedge\beta^{n-m}$, defined as*

$$\left[dd^cu\wedge dd^cv_1\wedge...\wedge dd^cv_{m-1}\wedge\beta^{n-m}\right](\omega)=\int u\, dd^cv_1\wedge...\wedge dd^cv_{m-1}\wedge\beta^{n-m}\wedge dd^c\omega,\ \omega\in F^{0,0} \qquad (5)$$

*is positive.*

The set of $m-sh$ in $D$ functions we denote by $sh_m(D)$. We note following properties of $m-sh$ functions (for more details, see [5]):

*1) if $u,v\in sh_m$, then $au+bv\in sh_m$ for any $a,b\geq 0$ ,i.e. the class $sh_m(D)$ represents convex cone;*

*2) $psh=sh_n\subset...\subset sh_1=sh$;*

*3) if $\gamma(t)-$convex, increasing function of parameter $t\in\mathbb{R}$ and $u\in sh_m$, then $\gamma\circ u\in sh_m$;*

*4) limit of uniformly convergent or decreasing sequence of $m-sh$ functions is $m-sh$;*

*5) maximum of finite numbers of $m-sh$ functions is $m-sh$ function ;*

*For arbitrary locally uniformly bounded family $\{u_\theta\}\subset sh_m$ the regularization $u^*(z)$ of supreme $u(z)=\sup_\theta u_\theta(z)$ also $m-sh$ function. Since $sh_m\subset sh$, then the set $\{u(z)<u^*(z)\}$ is polar in $\mathbb{C}^n\approx\mathbb{R}^{2n}$. Particularly, it has Lebesgue measure zero. Just as for locally uniformly bounded sequence $\{u_j\}\subset sh_m$ the regularization $u^*(z)$ of $u(z)=\overline{\lim_{j\to\infty}}u_j(z)$ also $m-sh$ function, at that the set $\{u(z)<u^*(z)\}$ is polar;*

*6) if $u\in sh_m$, then for any complex hyperplane $\mathrm{P}\subset\mathbb{C}^n$ the restriction $u|_{\mathrm{P}}\in sh_{m-1}$.*

**2.** One of the main problem of construction of potential theory in the class $sh_m(D)$ is to define operator $(dd^cu)^m\wedge\beta^{n-m}$ and introduction of capacity of condenser. We solve this problem on following scheme, which proposed by first author in alternative construction of pluripotential theory (see. [1,2]):

1) definition of operator $(dd^cu)^m\wedge\beta^{n-m}$ in class $sh_m(D)\cap C(D)$ (p.2);

2) $m-$polar set , $\mathscr{P}$*–measure* and their properties (p3,4);

3) definition of $m-$capacity $C_m(E,D)$ using just $sh_m(D)\cap C(D)$ (p.5);

4) proof of potential properties of $m-sh$ functions ( quasicontinuity, comparison principles, ets); definition of operator $(dd^c u)^m \wedge \beta^{n-m}$ in class $sh_m(D)\cap L^\infty_{loc}(D)$ and convergence $(dd^c u_j)^m \wedge \beta^{n-m} \mapsto (dd^c u)^m \wedge \beta^{n-m}$ for $u_j \downarrow u$ (p.6).

Let $1\le m\le n$ and $u_1,...,u_m \in sh_m(D)\cap C(D)$ . Then recurrence relation

$$\left[dd^c u_1 \wedge ... \wedge dd^c u_k \wedge \beta^{n-m}\right](\omega) = \int u_k dd^c u_1 \wedge ... \wedge dd^c u_{k-1} \wedge \beta^{n-m} \wedge dd^c\omega ,$$

$$\omega \in F^{m-k,m-k}, k = 1,...,m , \qquad (6)$$

defines positive current of bi-degree $(n-m+k, n-m+k)$ , at that for standard approximation $u_{ij} \downarrow u_i$ , $i = 1,2,...,k,\ j\to\infty$ we have convergence of currents

$dd^c u_{1j} \wedge ... \wedge dd^c u_{kj} \wedge \beta^{n-m} \mapsto dd^c u_1 \wedge ... \wedge dd^c u_k \wedge \beta^{n-m}$ (see [5]).

We note also, along with $dd^c u_1 \wedge ... \wedge dd^c u_k \wedge \beta^{n-m}$ in class $sh_m(D)\cap C(D)$ , just as defined a current $du_1 \wedge d^c u_1 \wedge dd^c u_2 \wedge ... \wedge dd^c u_k \wedge \beta^{n-m}$. It is easy to prove , that

$du_{1j} \wedge d^c u_{1j} \wedge dd^c u_{2j} \wedge ... \wedge dd^c u_{kj} \wedge \beta^{n-m} \mapsto du_1 \wedge d^c u_1 \wedge dd^c u_2 \wedge ... \wedge dd^c u_k \wedge \beta^{n-m}$ at $j\to\infty$.

Next integral estimation is very helpful in uniform estimations of $(dd^c u)^k \wedge \beta^{n-k}$ for the family of locally bounded $m-sh$ functions.

**Theorem 1.** *If* $u_1,u_2,...,u_k \in sh_m(B)\cap C(B)$, *where* $B=\{|z|<1\}$ *is ball and* $1\le k\le m$, *then for any* $r<1$

$$\int_0^r dt \int_{|z|^2\le t} dd^c u_1 \wedge dd^c u_2 \wedge ... \wedge dd^c u_k \wedge \beta^{n-k} \le (M-m) \int_{|z|^2\le r} dd^c u_2 \wedge ... \wedge dd^c u_k \wedge \beta^{n-k+1} ,$$

*where* $M = \sup_B u_1(z),\ m = \inf_B u_1(z)$ .

Proof of the theorem 1 is identical to proof of corresponding estimation for $psh$ functions [1].

**Corollary.** *In class of functions* $L_M = \{u \in psh(D)\cap C(D) : |u| \le M\}$ *the family of positive currents*

$\{dd^c u_1 \wedge ... \wedge dd^c u_k \wedge \beta^{n-m}\}, \{du_1 \wedge d^c u_1 \wedge dd^c u_2 \wedge ... \wedge dd^c u_k \wedge \beta^{n-m}\},\ u_1,...,u_k \in L_M, (1\le k\le m)$,

*weakly bounded.*

**3.** $m-$**polar sets.** By analogy of polar and pluripolar sets $m-$polar set will define as singular sets of $m-sh$ functions.

**Definition 3.** *A set* $E \subset D \subset \mathbb{C}^n$ *called* $m$*–polar in* $D$*, if there is a function* $u(z) \in sh_m(D)$, $u(z) \not\equiv -\infty$*, such that* $u\big|_E = -\infty$.

From $psh(D) \subset sh_m(D) \subset sh(D)$ follows that any pluripolar set is $m$–polar and in one's turn any $m$–polar set is polar. In particular, Hausdorff measure $H_{2n-2+\varepsilon}(E) = 0 \ \forall \varepsilon > 0$.

Now, we formulate a several theorems, which are identical to corresponding theorems for pluripolar sets, so we give these theorems without proofs.

**Theorem 2.** *Countable union of* $m$*–polar sets are* $m$*–polar, i.e. if* $E_j \subset D$ *are* $m$*–polar, then* $E = \bigcup_{j=1}^{\infty} E_j$ *is also* $m$*–polar.*

A domain $D \subset \mathbb{C}^n$ called *$m$–convex,* if there is a function $\rho(z) \in sh_m(D)$ such that $\lim_{z \to \partial D} \rho(z) = +\infty$; the domain D called *$m$–regular* if there is a function $\rho(z) \in sh_m(D)$, $\rho(z) < 0$ : $\lim_{z \to \partial D} \rho(z) = 0$.

**Theorem 3.** Let $D \subset \mathbb{C}^n$ *to be a* $m$*– convex domain and a subset* $E \subset D$ *is such, that for any compact domain* $G \subset\subset D$ *a set* $E \cap G$ *is* $m$*–polar in* $G$*. Then* $E$ *is* $m$*–polar in* $D$*. Moreover, if* $D$ *is* $m$*–regular, then there is a function* $u(z) \in sh_m(D)$, $u\big|_D < 0$, $u \not\equiv -\infty$*, such that* $u\big|_E \equiv -\infty$.

The theorem 3 is a preliminary results and we use it for proving more general result: local $m$–polar set is global (in $\mathbb{C}^n$) $m$–polar set.

**4. $\mathscr{P}$- measure.** Let $E \subset D$ is some subset of domain $D \subset \mathbb{C}^n$ and $1 \le m \le n$.For simplicity, below we suppose that $D$ is strong $m$–convex, i.e. $D = \{\rho(z) < 0\}$, where $\rho(z)$ is continuous and $m-sh$ function in some neighborhood $G \supset \overline{D}$.

**Definition 4.** *We consider class of function*

$$\mathcal{U} = \mathcal{U}(E,D) = \{u(z) \in sh_m(D) : u\big|_D \le 0,\ u\big|_E \le -1\} \text{ and put } \omega(z,E,D) = \sup_{u \in \mathcal{U}} u(z).$$

*Then the regularization* $\omega^*(z,E,D)$ *is called to be* $\mathscr{P}$*–measure (*$m$*–subharmonic measure) of the* $E \subset D$.

From property 7 of $m-sh$ function follows that $\omega^*(z,E,D) \in sh_m(D)$. $\mathscr{P}$–measure has following simple properties.

1) (Monotony) if $E_1 \subset E_2$*, then* $\omega^*(z,E_1,D) \ge \omega^*(z,E_2,D)$*, if* $E \subset D_1 \subset D_2$*, then* $\omega^*(z,E,D_1) \ge \omega^*(z,E,D_2)$*;*

2) *if* $U \subset D$ *–open set,* $U = \bigcup_{j=1}^{\infty} K_j$, *where* $K_j \subset \overset{\circ}{K}_{j+1}$, *then* $\omega^*(z, K_j, D) \downarrow \omega(z, U, D)$*;*

3) *if* $E \subset D$ *an arbitrary set, then there are a decreasing sequence of open sets* $U_j \supset E,\ U_j \supset U_{j+1}\ (j = 1, 2, ...)$, *such that* $\omega^*(z, E, D) = [\lim_{j\to\infty} \omega(z, U_j, D)]^*$*;*

4) $\mathscr{P}$ *–measure* $\omega^*(z, E, D)$ *is either nowhere zero or identically zero.* $\omega^*(z, E, D) \equiv 0$ *if and only if, when* $E$ *-is* $m-$*polar in* $D$.

**Definition 5.** *A point* $z^0 \in K$ *called* $m-$*regular point of compact* $K$ *(relatively* $D$*), if* $\omega^*(z^0, K, D) = -1$. *Compact* $K \subset D$ *called* $m-$*regular compact, if each point* $z^0$ *of* $K$ *is* $m-$*regular.*

Regular compact of classical potential theory are $m-$regular and $m-$regular compacts are pluriregular. It follows, that for any pair $K \subset U$, where $K-$compact and $U-$ open set, there is a $m-$regular compact $E$, such that $K \subset E \subset U$.

5) *if compact* $K \subset D$ *is* $m-$*regular, then* $\mathscr{P}-$ *measure* $\omega^*(z, K, D) \equiv \omega(z, K, D)$ *and is continuous function in* $D$. *Moreover, for* $m-$*regular compact* $\mathscr{P}-$ *measure* $\omega(z, K, D)$ *is maximal in* $D \setminus K$, $\left(dd^c\omega^*(z, K, D)\right)^m \wedge \beta^{n-m} = 0$*;*

**5. Condenser capacity. Definition 6.** *Let* $K \subset D \subset \mathbb{C}^n$. *Then a value*

$$C(K) = C(K, D) = \inf\left\{\int_D \left(dd^c u\right)^m \wedge \beta^{n-m} : u \in sh_m(D) \cap C(D),\ u|_K \le -1,\ \varliminf_{z\to\partial D} u(z) \ge 0\right\} \quad (7)$$

*is called capacity (*$m-$*capacity) of condenser* $(K, D)$.

The capacity has following properties:

1) *for* $m-$*regular compact* $K \subset D$ $\inf$ *in (7) reaches on* $\mathscr{P}$ *–measure, i.e.*

$$C(K) = \int_K \left(dd^c\omega^*(z, K, D)\right)^m \wedge \beta^{n-m};$$

By standard way we define exterior capacity assuming

$$C^*(E) = \inf\{C(U) : U \supset E - open\},$$

where capacity of open set

$$C(U) = \sup\{C(K) : K \subset U\} = \sup\{C(K) : K \subset U,\ K - regular\}.$$

2) *For any compact* $K \subset D$

$$C(K) = C^*(K) = \inf\{C(U) : U \supset K - open\} = \inf\{C(E) : E \supset K, E - regular\};$$

3) *if* $U \subset D$ *open set, then*

$$C(U) = \sup\{\int_U \left(dd^c u\right)^m \wedge \beta^{n-m} : u \in sh_m(D) \cap C(D),\ -1 \le u < 0\} =$$

$$= \sup\{\int_U \left(dd^c u\right)^m \wedge \beta^{n-m} : u \in sh_m(D) \cap C^\infty(D),\ -1 \le u < 0\}. \qquad (8)$$

Second supreme in (8) is useful, so, as integrand function is ordinary (regular).

4) *exterior capacity* $C^*(E)$ *monotonic, i.e. if* $E_1 \subset E_2$ *,then* $C^*(E_1) \le C^*(E_2)$*; it is countable -subadditive, i.e.* $C^*\left(\bigcup_j E_j\right) \le \sum_j C^*(E_j)$ *;*

5) *if* $E \subset D \subset G$, *then* $C^*(E,D) \ge C^*(E,G)$*;*

6) *for any increasing sequence of open sets* $U_j \subset U_{j+1}$ *holds* $C\left(\bigcup_j U_j\right) = \lim_{j\to\infty} C(U_j)$*;*

7) *exterior capacity of condenser* $C^*\left(E,D\right) = 0$ *if and only if, when* $E$ *is* $m-$*polar in* $D$*;*

**6.** Above introduced $\mathscr{P}$ – measure, condenser capacity and formulated their properties allow us to prove a several fundamental theorem of potential theory.

**Theorem 4.** *If a set* $E \subset \mathbb{C}^n$ *is locally* $m-$*polar, i.e. if for each point* $z^0 \in E$ *there is neighborhood* $B = B\left(z^0, r_{z^0}\right)$ *and a* $m-sh$ *in it function* $u(z) \not\equiv -\infty$, *such that* $u\big|_{E\cap B} \equiv -\infty$, *then* $E$ *is global* $m-$*polar* $\mathbb{C}^n$.

Theorem 4 for pluripolar sets using approximation of locally pluripolar sets with algebraic was proved by Josefsson [10]. In a work [1] proposed simple proof, based on condenser capacity, which passes also for Stein manifold. We give proof of theorem 4 using following chains.

÷ Fix a point $z^0 \in E$. Then there is a ball $B_r = B(z^0, r)$ such that $E \cap B_r$ is $m-$polar in $B_r$;

÷ $C^*(E \cap B_r, B_r) = 0$ (property 7 p.5) ;

÷ $C^*(E \cap B_r, B_R) = 0 \ \ \forall R > r$ (property 5 p.5) ;

÷ $E \cap B_r$ $m-$ polar in $B_R$ , $\forall R > r$ (property 7 p.5) ;

÷ $E \cap B_r$ $m-$polar in $C^n$, i.e. $\exists\ u_{z^0}(z) \in psh(C^n),\ \ u_{z^0} \not\equiv -\infty,\ \ u_{z^0}|_{E\cap B_r} \equiv -\infty$ (theorem 3) ;

÷ There are countable sets of such balls $B(z^j, r_j)$ covering $E$: $E \subset \bigcup_j B(z^j, r_j)$ and consequently, $E$ is $m-$polar in $C^n$. ⊳

Well-known C-property of N. N. Luzin confirms that any measurable function is continuous almost everywhere by Lebesgue measure. For $m-sh$ function one have continuity (quasicontinuity) almost everywhere by capacity (analogue of Cartan's theorem).

**Theorem 5.** *$m-$subharmonic function is continuous by capacity everywhere , i.e. if $u \in sh_m(D)$ , then for any $\varepsilon > 0$ there is an open set $U \subset D$ such that $C(U,D) < \varepsilon$ and u continuous in $D \setminus U$.*

Using this theorem we can proof next fundamental theorem of potential theory.

**Theorem 6.** *Let $1 \le m \le n$ and $u_1,...,u_m \in sh_m(D) \cap L^\infty_{loc}(D)$ . Then*

*1) the recurrence relation*

$$\left[dd^c u_1 \wedge ... \wedge dd^c u_k \wedge \beta^{n-m}\right](\omega) = \int u_k dd^c u_1 \wedge ... \wedge dd^c u_{k-1} \wedge \beta^{n-m} \wedge dd^c \omega , \quad \omega \in F^{m-k,m-k}(D), k = 1,...,m, \tag{9}$$

*defines positive current bi-degree $(n-m+k, n-m+k)$ ;*

2*) for a standard approximation $u_{ij} \downarrow u_i$ , $i = 0,1,...,m,\ j \to \infty$ we have convergence of currents $dd^c u_{1j} \wedge ... \wedge dd^c u_{kj} \wedge \beta^{n-m} \mapsto dd^c u_1 \wedge ... \wedge dd^c u_k \wedge \beta^{n-m}$ , $(k = 1,...,m)$;* (10)

In case, $u_1,...,u_m \in sh_m(D) \cap C(D)$ the proof easily follows from uniformly convergence $u_{ij} \downarrow u_i$ , $1 \le i \le m,\ j \to \infty$ , and in general case instead of continuity $u_i$ we have to use their quasicontinuity.

**Corallary.** *For any monotony decreasing sequence of $m-sh$ in $D$ functions $\{u_j(z)\}$ such that a limit $u(z) = \lim_{j\to\infty} u_j(z)$ locally bounded (from below) we have convergence of currents:*

*1)* $(dd^c u_j)^k \wedge \beta^{n-m} \mapsto (dd^c u)^k \wedge \beta^{n-m}$ ;

*2)* $u_j (dd^c u_j)^k \wedge \beta^{n-m} \mapsto u (dd^c u)^k \wedge \beta^{n-m}$ , $0 \le k \le m$.

**Remark.** In the paper [5] Z. Blocki proposed another method of definition of Hessian $(dd^c u)^m \wedge \beta^{n-m}$ in $sh_m(D)$. Let $\mathrm{D_m} \subset sh_m(D)$ class function $u \in sh_m(D)$ such that there is a Borel measure $\mu$, for which the current $(dd^c u_j)^m \wedge \beta^{n-m} \mapsto \mu\beta^n$ $\forall u_j \in sh_m(D) \cap C^2(D) : u_j \downarrow u$. He proved that $sh_m(D) \cap L^\infty_{loc}(D) \subset \mathrm{D_m}$. From mentioned above corollary follows that the measure $\mu$ for $u \subset sh_m(D) \cap L^{\infty}_{loc}(D)$ must be coincides with $(dd^c u)^m \wedge \beta^{n-m}$.

For $m-sh$ functions we also have comparison principle, which proved by Bedford and Taylor [6] for class of bounded $psh$ functions. We formulate it in following convenient form.

**Theorem 7.** (see also [5]). *If* $u,v \in sh_m(D) \cap L^{\infty}_{loc}(D)$ *and a set* $F = \{z \in D : u(z) < v(z)\} \subset\subset D$ *, then*

$$\int_F (dd^c u)^m \wedge \beta^{n-m} \geq \int_F (dd^c v)^m \wedge \beta^{n-m} \quad .$$

Geometrically theorem 7 means that in class $sh_m(D) \cap L^{\infty}_{loc}(D)$ operator $(dd^c v)^m \wedge \beta^{n-m}$ responsible for domination property. In particular, if $(dd^c v)^m \wedge \beta^{n-m} = 0$, then $v$ is maximal function.

**Theorem 8.** *For any compact $K \subset D$ its $\mathscr{P}$ – measure $\omega^*(z,K,D)$ satisfies in $D \setminus K$ the equation* $(dd^c \omega^*)^m \wedge \beta^{n-m} = 0$ .

We note that in p.4 (property 5) such fact reduced for $m-$regular compact $K \subset D$.

**Theorem 9.** The set $I_K$ *of irregular points of $K$ has zero capacity:* $C(I_K) = 0$*, i.e. $I_K$ is $m-$polar set.*

Next theorem has a connection with theorem 8 and gives positive answer to the second problem of Lelong for $m-sh$ functions.

**Theorem 10.** *Let $\{u_j\}$ is a increasing sequence of $m-sh$ functions such that* $u(z) = \lim_{j\to\infty} u_j(z)$ *is locally bounded from above. Then the set* $\sigma = \{u(z) < u^*(z)\}$ *is $m-$polar, where $u^*$ is regularization of $u$.*

A.S.Sadullaev. National University of Uzbekistan named after M.Ulugbek, *sadullaev@mail.ru*

B.I.Abdullaev. Urgench State University named after Al-Khorezmi, *abakhrom1968@mail.ru*